\documentclass[hidelinks]{ifacconf} %[letter, conference]
\usepackage{amsmath,amssymb,amsfonts}
\makeatletter
\let\old@ssect\@ssect % Store how ifacconf defines \@ssect
\makeatother

\usepackage{natbib}
\usepackage{hyperref}

\makeatletter
\def\@ssect#1#2#3#4#5#6{%
  \NR@gettitle{#6}% Insert key \nameref title grab
  \old@ssect{#1}{#2}{#3}{#4}{#5}{#6}% Restore ifacconf's \@ssect
}

\usepackage{parskip} %noident everywhere
\begingroup
    \makeatletter
    \@for\theoremstyle:=definition,remark,plain\do{%
        \expandafter\g@addto@macro\csname th@\theoremstyle\endcsname{%
            \addtolength\thm@preskip\parskip
            }%
        }
\endgroup
\usepackage{todonotes}

\usepackage[absolute]{textpos}
\usepackage{balance}
\usepackage{graphicx}
\usepackage{overpic}
\usepackage{rotating}
\usepackage{mathtools}
\usepackage{subcaption}
\usepackage{algorithm}
\usepackage{algpseudocode}

\usepackage{mymath}
\usepackage{nth}

\usepackage{scalerel}
\usepackage{tabularx}

\makeatletter
\renewcommand*\env@matrix[1][\arraystretch]{%
  \edef\arraystretch{#1}%
  \hskip -\arraycolsep
  \let\@ifnextchar\new@ifnextchar
  \array{*\c@MaxMatrixCols c}}
\makeatother

\usepackage{graphicx}
\usepackage{textcomp}
\usepackage{xcolor}
\def\BibTeX{{\rm B\kern-.05em{\sc i\kern-.025em b}\kern-.08em
    T\kern-.1667em\lower.7ex\hbox{E}\kern-.125emX}}

\definecolor{oxforduniversityblue}{RGB}{0,33,71}
\definecolor{oxforduniversityred}{RGB}{204, 41, 0}

\newcommand{\algorithmicbreak}{\textbf{break}}
\newcommand{\Break}{\State \algorithmicbreak}

\begin{document}
\begin{frontmatter}

\title{Fast Gradient Method for Model Predictive Control with Input Rate and Amplitude Constraints\thanksref{footnoteinfo}}

\thanks[footnoteinfo]{This research is supported by the Engineering and Physical Sciences Research Council (EPSRC) with a Diamond CASE studentship.}

\author[First]{Idris Kempf, Paul Goulart \& Stephen Duncan}
%\author[First]{ Paul Goulart}
%\author[First]{ Stephen Duncan}

\address[First]{Department of Engineering Science, University of Oxford, UK\\
   (e-mail: \{idris.kempf, paul.goulart, stephen.duncan\}@eng.ox.ac.uk)}
%\address[Second]{Colorado State University,
%   Fort Collins, CO 80523 USA (e-mail: author@lamar. colostate.edu)}
%\address[Third]{Electrical Engineering Department,
%   Seoul National University, Seoul, Korea, (e-mail: author@snu.ac.kr)}

\begin{abstract}                % Abstract of not more than 250 words.
This paper is concerned with the computing efficiency of model predictive control (MPC) problems for dynamical systems with both rate and amplitude constraints on the inputs. Instead of augmenting the decision variables of the underlying finite-horizon optimal control problem to accommodate the input rate constraints, we propose to solve this problem using the fast gradient method (FGM), where the projection step is solved using Dykstra's algorithm. We show that, relative to the Alternating Direction of Method Multipliers (ADMM), this approach greatly reduces the computation time while halving the memory usage. Our algorithm is implemented in C and its performance demonstrated using several examples.
\end{abstract}

\begin{keyword}
Model Predictive Control (MPC), Fast Gradient Method (FGM), Dykstra's Method, Alternating Direction of Multipliers Method (ADMM), Projection, Rate Constraints
\end{keyword}

\end{frontmatter}

%%%%%%%%%%%%%%%%%%%%%%%%%%%%%%%%%%%%%%%%%%%%%%%%%%%%%%%%%%%%%%%
% SECTION 1: Introduction, Paper Structure, Notation
%%%%%%%%%%%%%%%%%%%%%%%%%%%%%%%%%%%%%%%%%%%%%%%%%%%%%%%%%%%%%%%
\section{Introduction}
 Despite its advantages for constraint handling and feedforward disturbance modeling, the applicability of model predictive control is limited by the requirement to solve optimization problems in real-time to compute the control law.  While optimization problems are easily solved on a standard computer in a simulation environment, it is considerably more challenging to solve such problems on the embedded systems often encountered in industrial applications. Cost factors often outweigh the need for powerful hardware to implement an MPC scheme using standard optimization program solvers. Embedded systems employed in industrial applications, e.g.\ in the aviation or automotive industry, therefore feature considerably lower computational power than a standard desktop computer and much tighter memory restrictions.  Hardware platforms in the aviation and automotive industry commonly run at clock frequencies of a few hundred megahertz, while offering only a few megabytes of memory, whereas standard computers run at gigahertz rates and can provide many gigabytes of memory storage. A key to employing MPC in industrial applications is therefore an algorithm that makes maximum use of available computing resources while keeping the memory usage at a minimum.

 A common theme in deterministic MPC formulations encountered in practice is that, in order to avoid infeasible optimization problems, constraints are most often imposed only on the inputs rather than on the states of the dynamical system. The reason is that any disturbance acting on the system might drive the states outside of the feasible constraint set, thereby resulting in an infeasible optimization problem. Moreover, it is also the case that two different kinds of constraints are imposed on the input: rate and amplitude constraints. This should not be surprising as most system inputs drive some kind of actuator which naturally possesses some performance limits.

 Algorithms that aim to solve constrained optimal control problems using first-order techniques typically require a projection\footnote{Throughout the paper, it is assumed that an Euclidean projection is used.} onto the constraint set. In the absence of rate constraints, this projection is usually straightforward and can be computed using a closed-form formula \citep[Chapter 28.3]{BAUSCHKEBOOK}, e.g.\ projection onto a box-shaped set of upper and lower actuator limits. However, if input rate constraints (often referred to as \emph{slew rate}) are included then the projection is more complicated, and we are not aware of any closed form solution to the corresponding projection problem (see also \citet{BAUSCHKESWISS}).  Most projection-based optimization methods therefore employ an augmented problem form, e.g.\ one that includes additional state variables. The projection problem then reverts to a closed-form formula. While this approach is versatile in the sense that it can cope with most reasonable sets encountered in practice, the augmentation of decision variables curtails the computation speed while increasing the memory usage. It also introduces additional equality constraints to the problem, leading to difficulties in applying methods such as the fast gradient method. The question arises whether it is actually necessary to augment the decision variables in the particular case of constraints arising from  input rate and amplitude constraints.

 This paper suggests an approach that does \emph{not} require augmenting the decision variables of the optimization problem.  By employing a closed-form solution for a 2-dimensional rate and amplitude constraint set in combination with Dykstra's algorithm \citep{DYKSTRA}, we show that the projection of a vector of arbitrary finite dimension onto the space of rate-constraint signals can be found iteratively, obviating the need for additional state or other problem variables. Our projection algorithm is then embedded in a fast gradient method \citep{OPTIMNESTEROV}. It turns out that the complexity introduced by our projection algorithm is considerably less than the complexity introduced by the augmentation of decision variables.

 This paper is organized as follows. In Section \ref{chap:problem}, we introduce both the linear MPC problem and the fast gradient method - an algorithm which is particularly suitable for solving it.  Section \ref{chap:projection} formally defines the input rate and amplitude constraint set before presenting Dykstra's method and its application to the former set. In Section \ref{chap:mpcfgm}, the fast gradient method is combined with Dykstra's method. The performance of our algorithm is compared with the Alternating Direction of Method of Multipliers  (ADMM) for MPC problems with input rate and amplitude constraints. Both algorithms are implemented in C and tested on several examples.

\textit{Notation} Let $\otimes$ denote the Kronecker product and $\oplus$ the direct sum (i.e.\ the block diagonal concatenation) of two matrices, respectively. Let $I_n$ represent the identity matrix in $\R^{n\times n}$. Let $\twonorm{\cdot}$ and $\infnorm{\cdot}$ denote the two- and infinity-norm, respectively.

%%%%%%%%%%%%%%%%%%%%%%%%%%%%%%%%%%%%%%%%%%%%%%%%%%%%%%%%%%%%%%%
% PROBLEM STATEMENT
%%%%%%%%%%%%%%%%%%%%%%%%%%%%%%%%%%%%%%%%%%%%%%%%%%%%%%%%%%%%%%%
\section{Problem Statement}\label{chap:problem}
\subsection{Model Predictive Control}\label{chap:mpc}
Given a discrete-time linear dynamical system and an initial condition $x(t)$ at time $t$, a standard MPC scheme computes a control law by predicting the future evolution of the system and minimizing a quadratic objective function over some planning horizon $T$. This can be achieved via repeated solution of the following quadratic program (QP):
\begin{subequations}
\begin{align}
\min &\sum_{k=0}^{T-1} x_k^\Tr Q x_k+u_k^\Tr R u_k + x_T^\Tr P x_T,\label{eq:mpcA}\\
s.t.\,\,\, &x_{k+1} = Ax_k+Bu_k, \quad\, x_0=x(t),\label{eq:mpcB}\\
     %&(x_1,\dots,x_{T-1})\in\mathcal{X}_k,\label{eq:mpcD}\\
     &(u_1,\dots,u_{T-1})\in\mathcal{U},\label{eq:mpcC}
\end{align}\label{eq:mpc}\end{subequations}
for $k=0,\dots,T-1$, returning at each step the optimal first input stage $u^*_0=u(t)$ as a control law. The inputs $u_k\in\R^{n_u}$ are constrained to the closed convex set $\mathcal{U}\subseteq\R^{n_u\times n_u}$. It is assumed that no constraints are imposed on the states $x_k\in\R^{n_x}$. The terminal cost matrix $P\!=\!P^\Tr\!\succ\! 0$ is obtained from the discrete-time algebraic Riccati equation (DARE) associated to the unconstrained infinite horizon regulator problem. The QP \eqref{eq:mpc} has a unique solution if $R\!\succ\!0$, $Q\!\succeq\!0$ and the pairs $(A,B)$ and $(A,Q^{\frac{1}{2}})$ are controllable and observable, respectively \citep[Chapter 12]{MPCBOOK}.
%\PG{We do not include any terminal state constraint to ensure stability as in \cite{mayne2000constrained}, since such constraints are typically neglected in industrial applications.   It is always possible to ensure stability \emph{without} such a constraint if $T$ is sufficiently large \citep{jadbabaie2005stability}.}

 By eliminating the state variables $x\eqdef (x_1, \dots, x_T)$ and defining $u \eqdef (u_0,\dots,u_{T-1})^\Tr$, \eqref{eq:mpc} can be reformulated in its condensed form as
\begin{subequations}\label{eq:mpcqp}
\begin{align}
\min \,\,\,&\frac{1}{2}u^\Tr J u +q^\Tr u,\label{eq:mpcqpA}\\
s.t.\,\,\, &u\in\mathcal{U},\label{eq:mpcqpB}
\end{align}\end{subequations}
where $J$ and $q$ are defined as
\begin{subequations}\label{eq:Jq}
\begin{align}
J &\eqdef G^\Tr \left( (I_T\otimes Q)\oplus P \right)G+(I_T\otimes R),\label{eq:JqA}\\
q &\eqdef G^\Tr H x_0,\label{eq:JqB}
\end{align}\end{subequations}
and $G$ and $H$ arise from elimination of the equality constraints in \eqref{eq:mpcB}, i.e.\ from writing \eqref{eq:mpcB} as $x=Gu+Hx_0$. Note that $J\!=\!J^\Tr\! \succ\!0$ by the assumptions of the previous paragraph.

\subsection{Fast Gradient Method}\label{chap:fgm}
The fast gradient method (FGM) belongs to the family of first-order methods that seek solutions of convex optimization problems using only the first derivative of the objective function \eqref{eq:mpcqpA}. In this paper we will use the \textit{constant step scheme II} \citep[Chapter 2.2]{OPTIMNESTEROV} which is known to have an \textit{optimal} convergence rate \citep[Thm. 2.2.2]{OPTIMNESTEROV}, but our results will hold generally for any first-order optimisation scheme employing a projection. A formulation of this algorithm is presented in \citep{FPGAMPC} and repeated in Algorithm \ref{alg:fgm1}. The fixed step size $\beta = (\sqrt{\lambda_{max}}-\sqrt{\lambda_{min}})/(\sqrt{\lambda_{max}}+\sqrt{\lambda_{min}})$ is based on the minimum and maximum curvature of the convex objective function, implying that the objective must be strongly convex. For \eqref{eq:mpcqp}, finding these values amounts to computing the minimum and maximum eigenvalues $\lambda_{min}$ and $\lambda_{max}$ of $J\succ 0$, respectively.
%\begin{align}
%\beta = \frac{\sqrt{\lambda_{max}}-\sqrt{\lambda_{min}}}{\sqrt{\lambda_{max}}+\sqrt{\lambda_{min}}},
%\end{align}

 In order to reduce the computation effort and complexity, we will not use any termination criterion and instead run the algorithm for a fixed number of iterations $I_{max}$. Based on the convergence rate results for the fast gradient method, a maximum number of iterations $I_{max}$ can be derived that guarantees a certain level of suboptimality for all initial states $x_0=x(t)$ within a bounded set \citep{CERTIFMPC}.

 In order to apply Algorithm \ref{alg:fgm1} to our optimal control problem \eqref{eq:mpcqp}, the projection operator $\mathcal{P}_{\mathcal{U}}$ for the set $\mathcal{U}$ must be known. For a constraint set combining both input rate and amplitude constraints, this projection is not straightforward and will be addressed in Section~\ref{chap:projection}.
\begin{algorithm}[H]
 \caption{Fast Gradient Method applied to MPC problem \eqref{eq:mpc}}\label{alg:fgm1}
 \begin{algorithmic}[1]
    \Require $x_0=x(t)$
    \Ensure $u(t)=u_{\scaleto{I}{4pt}_{max}}$
	\State Set $q=G^\Tr H x_0$ and $y_1=u_1=0$
	\For{$i = 1$ to $I_{max}$}
  	  \State $t_i = (I - J/\lambda_{max})y_i - q / \lambda_{max}$
  	  \State $u_{i+1} = \mathcal{P}_{\mathcal{U}}(t_i)$
  	  \State $y_{i+1} = (1+\beta)u_{i+1} - \beta u_i$
	\EndFor
 \end{algorithmic}
\end{algorithm}

\subsection{Alternating Direction of Multipliers Method}\label{chap:admm}
The Alternating Direction of Method of Multipliers (ADMM) belongs to the class of augmented Lagrangian methods and is, like FGM, a first-order method. ADMM algorithms are based on repeatedly minimizing the augmented Lagrange function w.r.t.\ the primal variables and maximizing the same function w.r.t.\ to the dual variables \citep[Chapter 5]{BOYDCONVEX}. Assuming that the input constraint set \eqref{eq:mpcqpB} can be represented as a polyhedron, i.e.\
\begin{align}
\mathcal{U}=\lbrace u \in \R^{T-1} \vert \ubar{v} \leq K u \leq \bar{v} \rbrace\label{eq:constraintspolyhedron}
\end{align}
the optimization problem \eqref{eq:mpcqp} can be reformulated as
\begin{subequations}\label{eq:mpcqpaug}
\begin{align}
\min \,\,\,&\frac{1}{2}u^\Tr J u +q^\Tr u,\label{eq:mpcqpaugA}\\
s.t.\,\,\, &K u - v = 0, \label{eq:mpcqpaugB}\\
& \ubar{v} \leq v \leq \bar{v},\label{eq:mpcqpaugC}
\end{align}\end{subequations}
where the constraint variables $v\in \R^{n_v}$  and equality constraints \eqref{eq:mpcqpaugB} were introduced. The augmented Lagrangian for \eqref{eq:mpcqpaug} can be written as
\begin{align}\label{eq:lang}
\begin{split}
L(u,v,\gamma) = &\frac{1}{2}u^\Tr J u + q^\Tr u+\frac{\rho}{2}\lVert Ku-v\rVert_2^2\\&+\gamma^\Tr(Ku-v) + \mathcal{I}_{[\ubar{v},\bar{v}]}(v),
\end{split}
\end{align}
where $\mathcal{I}_{[\ubar{v},\bar{v}]}$ is the indicator function for the set $\mathcal{V}=\set{v}{\ubar{v}\le v \le \bar{v}}$ and the penalty parameter $\rho>0$ and the dual variables $\gamma$ are associated with the constraint \eqref{eq:mpcqpaugB}. A standard ADMM scheme solves \eqref{eq:mpcqpaug} by repeatedly minimizing \eqref{eq:lang} w.r.t.\ $u$ and $v$ and updating the dual variables $\gamma$ using an approximate gradient ascent method. The algorithm is summarized in Algorithm \ref{alg:admm}, where the saturation function $\sat_{[\ubar{v},\bar{v}]}(\cdot)$ was used, which limits its argument to $\ubar{v}$ and $\bar{v}$. Reformulation \eqref{eq:mpcqpaug} simplifies the projection involved in the algorithm: Instead of projecting onto the polyhedron $\mathcal{U}$, which might be as hard as solving \eqref{eq:mpcqp}, the projection onto $\mathcal{V}$ is given by the saturation function. For the purposes of this paper, the critical distinction between our optimisation problem in the form \eqref{eq:mpcqp} (and its solution via FGM) and \eqref{eq:mpcqpaug} (and its solution via ADMM) is that the form \eqref{eq:mpcqp} has a positive definite $J$ and does not allow for equality constraints. The problem form \eqref{eq:mpcqpaug} makes neither restriction.

\begin{algorithm}[H]
 \caption{Alternating Direction of Multipliers Method applied to MPC problem \eqref{eq:mpc} with set \eqref{eq:constraintspolyhedron}}\label{alg:admm}
 \begin{algorithmic}[1]
 \Require{$x_0=x(t)$}
 \Ensure{$u(t)=u_{\scaleto{I}{4pt}_{max}}$}
  \State Set $q=G^\Tr H x_0$ and $\gamma_0 = 0$
  \For{$i = 1$ to $I_{max}$}
  	\State Solve for $u_i$:\\ $\left( J+\rho K^\Tr K \right) u_i = K^\Tr (\rho v_{i-1}-\gamma_{i-1}) -q$
  	\State $v_i = \sat_{[\ubar{v},\bar{v}]}\left\lbrace Ku_i+\rho^{-1}\gamma_{i-1}\right\rbrace$
  	\State $\gamma_i = \gamma_{i-1} +\rho(Ku_i-v_i)$
  \EndFor
 \end{algorithmic}
\end{algorithm}

%%%%%%%%%%%%%%%%%%%%%%%%%%%%%%%%%%%%%%%%%%%%%%%%%%%%%%%%%%%%%%%
% SECTION PROJECTION
%%%%%%%%%%%%%%%%%%%%%%%%%%%%%%%%%%%%%%%%%%%%%%%%%%%%%%%%%%%%%%%
\section{Input Constraint Projection Method}\label{chap:projection}
Given a nonempty closed convex set $\mathcal{X}\subseteq\R^{N\times N}$, the Euclidean projection $x^*$ of a point $x^\circ\in\R^{N}$ is defined as the minimizer of the following optimization problem:
\begin{align}\label{eq:projection}
x^* = \argmin_{x\in\mathcal{X}} \twonorm{x - x^\circ}.
\end{align}
By the assumptions on the set $\mathcal{X}$, the optimization problem \eqref{eq:projection} admits a unique solution \citep[Chapter 3.2]{BAUSCHKEBOOK}. In the following, we will write $x^*=\mathcal{P}_{\mathcal{X}}(x^\circ)$ as shorthand for projecting a point $x^\circ$ onto the set $\mathcal{X}$.

\subsection{Rate and Amplitude Constraint Set}\label{chap:constraints}
Given a maximum allowable input amplitude $a\!\in\!\R^{n_u}\!>\!0$ and rate $r\!\in\!\R^{n_u}\!>\!0$, define the amplitude constraint set as
$\mathcal{A}_k \eqdef \left\lbrace u \in \R^{n_u(T-1)}\,\vert\,\abs{u_k} \leq a\right\rbrace$ and the rate constraint set as $\mathcal{R}_k \eqdef \left\lbrace u \in \R^{n_u(T-1)}\,\vert\,\abs{u_k - u_{k-1}} \leq r\right\rbrace$ for $k=0,\dots,T-1$ and where the inequalities are applied element-wise. The set $\mathcal{R}_0$ includes the input $u_{-1}$ that is treated as a fixed constant stemming from the actual input of the system at time $t-1$.  We exclude the trivial case where $\mathcal{A}_k$ is entirely contained in $\mathcal{R}_k$ by assuming that $0\leq r \leq 2a$ for all $n_u$ elements. The input rate and amplitude constraint set for problem \eqref{eq:mpcqp} is obtained as the intersection of $\mathcal{A}=\mathcal{A}_0\cap\dots\cap\mathcal{A}_{T-1}$ and $\mathcal{R}=\mathcal{R}_0\cap\dots\cap\mathcal{R}_{T-1}$, i.e.\
\begin{gather}\label{eq:israc}
\begin{aligned}
\mathcal{U} \eqdef \bigl\lbrace  u \in \R^{n_u(T-1)}\,\bigl|\,&\abs{u_k} \leq a \,\,\forall\,\, k\in\mathcal{K},\\
&\abs{u_k - u_{k-1}} \leq r \,\,\forall\,\, k\in\mathcal{K}\bigr\rbrace.
\end{aligned}
\end{gather}
Because the constraints are not coupled among the elements of $u_k\in\R^{n_u}$, we will assume throughout that $n_u=1$ for clarity of exposition. However, all of our results apply in the case that $n_u>1$.

 While it is straightforward to obtain a closed-form formula for the projection onto $\mathcal{A}$ \citep[Chapter 28.3]{BAUSCHKEBOOK}, we know of no tractable closed-form solution for $\mathcal{P}_\mathcal{R}$ (see also \citet{BAUSCHKESWISS}) and hence also not for $\mathcal{P}_\mathcal{U}$. In order to see why this projection is difficult, define $N\eqdef T-1$ and consider the reformulation of $\mathcal{P}_\mathcal{R}(u^\circ)$ with $n_u=1$ as a dynamic programming problem \citep{BERTSEKAS}:
\begin{subequations}\label{eq:dpproblem}
\begin{align}
\min_{u,v}\,\,\, &\sum_{k=0}^{N} (u_k - u_k^\circ)^2,\\
s.t.\,\,\, & u_k = u_{k-1} + v_{k-1}, \quad\, u_{\sm 1}=u(t-1),\\
		   & -r \leq v_{k-1} \leq r, \quad k=0,\dots,N,
\end{align}\end{subequations}
where we have introduced the slope variable $v_{k-1}\eqdef u_k - u_{k-1}$. Adopting the notation from \citet{BERTSEKAS} and denoting the cost-to-go at step $k$ by $J_k(u_k)$, the solution to \eqref{eq:dpproblem} for the first two steps is obtained as:
\begin{itemize}
\item Step $k=N:$
\begin{align}
J_{N}(u_{N}) = (u_{N} - u_{N}^\circ)^2
\end{align}
\end{itemize}
\begin{itemize}
\item Step $k=N-1:$
\end{itemize}
\begin{subequations}
\begin{align}
&J_{N-1}(u_{N-1}) = \min_{u_{N-1}} \left\lbrace (u_{N-1} \!-\!u_{N-1}^\circ)^2\!+\!J_{N}(u_{N})\right\rbrace\\
&\implies v_{N-1}\begin{cases}
r\phantom{-u_{N}^\circ u_{N-1}}\,\text{ if } u_{N}^\circ - u_{N-1} \geq r\\
-r\phantom{u_{N}^\circ u_{N-1}}\,\text{ if } u_{N}^\circ - u_{N-1} \leq -r\\
u_{N}^\circ - u_{N-1}\,\text{ otherwise}
\end{cases}\label{eq:dpsol}
\end{align}\end{subequations}
At step $k=N-2$, each of the three cases of \eqref{eq:dpsol} induce three other conditions for $v_{N-2}$. Continuing this way, this amounts to $3^N=3^{T-1}$ conditions for obtaining a solution to \eqref{eq:dpproblem}. While it can be possible to derive a formula for small horizons, the solution via dynamic programming becomes intractable for larger horizons.

 Another approach to solve \eqref{eq:dpproblem} or directly the projection onto the rate and amplitude constraint set would be to define a multi-parametric program \citep[Chapter 2]{MPCBOOK} with parameters $u^\circ$ and $u_{\sm 1}$. The multi-parametric solution of \eqref{eq:dpproblem} results in a piecewise affine function (PWA) of the parameters $u_{\sm 1}$ and $u^\circ$, i.e.\ $n$ affine functions defined on $n$ disjoint sets. While an explicit solution could be computed using dedicated software, e.g.\ \citep{MPT3}, it is expected that the number of regions $n$ would be at least $3^{T-1}$ for projecting onto the set \eqref{eq:israc} with $n_u=1$.

\subsection{2-Dimensional Projection}\label{chap:2dim}
\begin{figure}
\begin{center}
\begin{overpic}[scale=1]{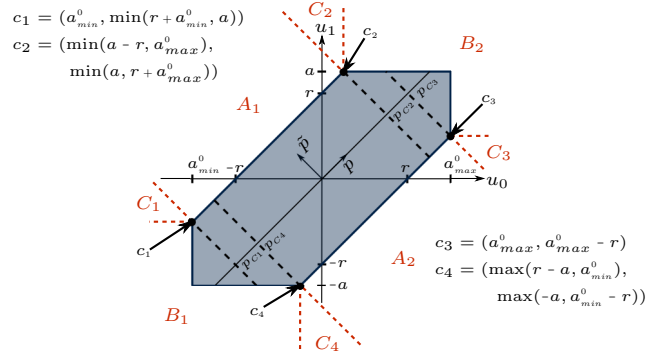}
		% Axes limits
		\put(52,17.5){\tiny{$\sm a$}}
		\put(52,24){\tiny{$\sm r$}}
		\put(45,80.5){\tiny{$a$}}
		\put(45,74){\tiny{$r$}}
		\put(11,53){\tiny{$a_{\scaleto{min}{2pt}}^{\scaleto{0}{2pt}}$}}
		\put(22,52){\tiny{$\sm r$}}
		\put(85.5,53){\tiny{$a_{\scaleto{max}{1.3pt}}^{\scaleto{0}{2pt}}$}}
		\put(73,52){\tiny{$r$}}
		% Coordinate axes
		\put(98,48){$\scriptstyle{u_0}$}
		\put(48,93){$\scriptstyle{u_1}$}
		\put(56,52){\rotatebox{45}{$\scriptstyle{p}$}}
		\put(42,58){\rotatebox{45}{$\scriptstyle{\tilde{p}}$}}
		% Diag axes limits
		\put(26.5,23.5){\rotatebox{45}{$\scriptscriptstyle{p_{\scaleto{C1}{2pt}}}$}}
		\put(33,29){\rotatebox{45}{$\scriptscriptstyle{p_{\scaleto{C4}{2pt}}}$}}
		\put(70.5,67){\rotatebox{45}{$\scriptscriptstyle{p_{\scaleto{C2}{2pt}}}$}}
		\put(77.5,74){\rotatebox{45}{$\scriptscriptstyle{p_{\scaleto{C3}{2pt}}}$}}
		% Regions
		\put(-4,41){\textcolor{oxforduniversityred}{$\scriptstyle{C_1}$}}
		\put(46,98){\textcolor{oxforduniversityred}{$\scriptstyle{C_2}$}}
		\put(98,56){\textcolor{oxforduniversityred}{$\scriptstyle{C_3}$}}
		\put(48,0){\textcolor{oxforduniversityred}{$\scriptstyle{C_4}$}}
		\put(4,8){\textcolor{oxforduniversityred}{$\scriptstyle{B_1}$}}
		\put(90,88){\textcolor{oxforduniversityred}{$\scriptstyle{B_2}$}}
		\put(25,70){\textcolor{oxforduniversityred}{$\scriptstyle{A_1}$}}
		\put(70,25){\textcolor{oxforduniversityred}{$\scriptstyle{A_2}$}}
		% Corners
		\put(-4,27.5){\tiny{$c_{\scaleto{1}{2pt}}$}}
		\put(62,92){\tiny{$c_{\scaleto{2}{2pt}}$}}
		\put(96.5,72.5){\tiny{$c_{\scaleto{3}{2pt}}$}}
		\put(29,10){\tiny{$c_{\scaleto{4}{2pt}}$}}
		% Corner Equations
		\put(-40,96){\tiny{$c_1 = (a_{\scaleto{min}{2pt}}^{\scaleto{0}{2pt}}, \min(r\,\scaleto{+}{3pt}\,a_{\scaleto{min}{2pt}}^{\scaleto{0}{2pt}}, a))$}}
		\put(-40,88){\tiny{$c_2 = (\min(a\sm r, a_{\scaleto{max}{2pt}}^{\scaleto{0}{2pt}}),$}}
		\put(-40,80){\tiny{$\phantom{c_2 = (} \min(a, r\,\scaleto{+}{3pt}\,a_{\scaleto{max}{2pt}}^{\scaleto{0}{2pt}}))$}}

		\put(83,30){\tiny{$c_3 = (a_{\scaleto{max}{2pt}}^{\scaleto{0}{2pt}}, a_{\scaleto{max}{2pt}}^{\scaleto{0}{2pt}} \sm r)$}}
		\put(83,22){\tiny{$c_4 = (\max(r \sm a, a_{\scaleto{min}{2pt}}^{\scaleto{0}{2pt}}),$}}
		\put(83,14){\tiny{$\phantom{C_4 = (}\max(\sm a, a_{\scaleto{min}{2pt}}^{\scaleto{0}{2pt}}\sm r))$}}
\end{overpic}
\end{center}\caption{Input rate and amplitude constraint set for the 2-dimensional projection as defined in \eqref{eq:setsplit0}. The regions $A_i, B_i$ and $C_i$ define different maps used to project onto the shaded set.}\label{fig:2dimset}
\end{figure}

Consider the input rate and amplitude constraint set \eqref{eq:israc} for $n_u=1$ and $T=2$, i.e.\ $u=(u_0, u_1)^\Tr$, which can be represented as
\begin{gather}\label{eq:setsplit0}
\begin{aligned}
\mathcal{U}_1 \eqdef \bigl\lbrace  u \in \R^{n_u(T\sm 1)}\,\bigl|\,&a^0_{min} \leq u_0 \leq a^0_{max},\\
&\abs{u_1} \leq a,\,\,\abs{u_1 - u_0} \leq r \bigr\rbrace,
\end{aligned}\end{gather}
where $a^0_{min}\eqdef \max(-a, -r+u_{\sm 1})$ and $a^0_{max}\eqdef \min(a, r+u_{\sm 1})$. The set is illustrated in Figure \ref{fig:2dimset}, where a coordinate system $(p,\tilde{p})$ rotated by $45^\circ$ and different regions $A_i, B_i$ and $C_i$ have been added. In addition, the corner points $c_i$ are defined and their projections onto the $p$-diagonal have been marked by $p_{c_i}$. When $r$ and $a$ are fixed, the shape of set \eqref{eq:israc} depends on parameter $u_{\sm 1}$. As $u_{\sm 1}$ changes, corner points $c_i$ are moved along the rate constraint diagonals until they eventually stop at an amplitude constraint boundary. In addition, regions $B_i$ might lose their horizontal facet and regions $A_i$ vanish. Assume that we wish to project a point $u^\circ=(u_0^\circ,u_1^\circ)^\Tr$. The projection consists of two steps: Firstly, determine to which region the point belongs. Secondly, apply the projection function of that particular region.  In other words, we identify explicitly the PWA solution to the problem of projection onto the set $\mathcal{U}_1$ following the general method of \citep[Chapter 2]{MPCBOOK}.
%Following conditions are obtained for the point location:
%\begin{subequations}\begin{align}
%(p_{c_1} \leq p_{u^\circ} \leq p_{c_2}) \,\&\, (\tilde{p}_{u^\circ} \geq r/\sqrt{2}) &\Longrightarrow A_1,\\
%(p_{c_4} \leq p_{u^\circ} \leq p_{c_3}) \,\&\, (\tilde{p}_{u^\circ} \leq r/\sqrt{2}) &\Longrightarrow A_2,\\
%(p_{u^\circ} \leq p_{c_1}) \,\&\, (u^\circ_0 \leq (c_1)_0) \,\&\, (u^\circ_1 \geq (c_1)_1) &\Longrightarrow C_1,\\
%(p_{u^\circ} \geq p_{c_2}) \,\&\, (u^\circ_0 \geq (c_2)_0) \,\&\, (u^\circ_1 \geq (c_2)_1) &\Longrightarrow C_2,\\
%(p_{u^\circ} \geq p_{c_3}) \,\&\, (u^\circ_0 \geq (c_3)_0) \,\&\, (u^\circ_1 \leq (c_3)_1) &\Longrightarrow C_3,\\
%(p_{u^\circ} \leq p_{c_4}) \,\&\, (u^\circ_0 \geq (c_4)_0) \,\&\, (u^\circ_1 \leq (c_4)_1) &\Longrightarrow C_4,
%\end{align}\label{eq:conditions}\end{subequations}
%where $(c_i)_0$ and $(c_i)_1$ denote the components of point $c_i$ in the $u_0$- and $u_1$-direction, respectively, $\&$ the Boolean addition and $p_{u^\circ}$ and $\tilde{p}_{u^\circ}$ the projection of $u^\circ$ onto the $p$- and $\tilde{p}$-diagonal, respectively. If none of the above conditions are met, the point lies in $B_1$, $B_2$ or inside the shaded area. Note that omitting the constraint $\abs{u_0 - u_{-1}}\leq r$ introduces a symmetry and considerably simplifies the point location and therefore the projection.

 After representing $u^\circ$ in the $(p,\tilde{p})$-basis, the point location problem can be easily solved using Figure \ref{fig:2dimset} and by comparing the $p$-component with the corner projections $p_{c_i}$ determined by $r$, $a$ and $u_{-1}$. If the $p$-component is beyond these limits, the $u_0$- or $u_1$-components can be compared with the corner points $c_i$ to complete the point location problem.

 If the point lies in one of the corner regions $C_i$, then it is mapped to the corresponding corner $c_i$. If the point lies in one of the diagonal regions $A_i$, then it can be represented in terms of the rotated basis $(p,\tilde{p})$, its $\tilde{p}$-component saturated to $\pm r/\sqrt{2}$ and rotated back. Finally, if the point lies in one of the box-regions $B_i$, then the formula for a box-projection can be applied \citep[Chapter 28.3]{BAUSCHKEBOOK} using limits $(a^0_{min},a^0_{max})$ and $(a^1_{min},a^1_{min})$ for the $u_0$- and $u_1$-direction, respectively, where $a^1_{min}= \max(-a, -r + a^0_{min})$ and $a^1_{max}= \min(a, r+a^0_{max})$. A complete C-language implementation of the 2-dimensional projection can be found in \citep{BRPGITHUB}.

\subsection{Dykstra's Algorithm}\label{chap:dykstra}
 Dykstra's algorithm \citep{DYKSTRA} was first published in 1983 as an extension to Von Neumann's \textit{Alternating Projections Method} \citep{NEUMANN}, which is suitable for finding a point lying in the intersection $\mathcal{U}=\mathcal{U}_1\cap\dots\cap\mathcal{U}_N$ of $N$ closed convex sets $\mathcal{U}_i$ by cyclically projecting onto the sets $\mathcal{U}_i$. While Von Neumann's algorithm only finds \textit{some} point in $\mathcal{U}$, Dykstra's algorithm determines the Euclidean projection $u^*=\mathcal{P}_{\mathcal{U}}(u^\circ)$ of $u^\circ$ onto $\mathcal{U}$. Both algorithms circumvent the potentially complicated projection $\mathcal{P}_{\mathcal{U}}$ by iteratively applying the (known) projections $\mathcal{P}_{\mathcal{U}_i}$. The method is summarized in Algorithm \ref{alg:dykstra} for the case that $\mathcal{U}=\mathcal{U}_1\cap\mathcal{U}_2$.
\begin{algorithm}[H]
 \caption{Dykstra's Algorithm for two sets}\label{alg:dykstra}
 \begin{algorithmic}[1]
  \Require{$v^\circ$}
  \Ensure{$\mathcal{P}_{\mathcal{U}}(v^\circ)=x_{i+1}$}
  \State{Set $x_0=v^\circ$, $\mu_0=0$ and $\gamma_0=0$}
  \For{$i = 1$ to $I_{max}$}
   \State{$y_i=\mathcal{P}_{\mathcal{U}_1}(x_i+\mu_i)$}
   \State{$\mu_{i+1}=\mu_i+x_i-y_i$}
   \State{$x_{i+1}=\mathcal{P}_{\mathcal{U}_2}(y_i+\gamma_i)$}
   \State{$\gamma_{i+1}=\gamma_i+y_i-x_{i+1}$}
   \If{$\infnorm{x_{i+1}-x_i} < \epsilon$}
  		\Break
  	\EndIf
  \EndFor
 \end{algorithmic}
\end{algorithm}
 What distinguishes Dykstra's algorithm from Von Neumann's is the choice of variables $\mu^k$ and $\gamma^k$ that track the residuals from projecting onto $\mathcal{U}_1$ and $\mathcal{U}_2$ \citep{DYKSTRAADMM}. It can be shown that Dykstra's algorithm is equivalent to  ADMM applied to problem \eqref{eq:projection} for the case that $\mathcal{U}=\mathcal{U}_1\cap\mathcal{U}_2$ \citep{DYKSTRAADMM}. In addition, it has been proved that Algorithm \ref{alg:dykstra} always converges to the Euclidean projection onto $\mathcal{U}$, provided that the sets $\mathcal{U}_1$ and $\mathcal{U}_2$ are closed convex sets and their intersection is nonempty \citep{DYKSTRA,BAUSCHKEBOOK}.

 Before applying Algorithm \ref{alg:dykstra} to the rate and amplitude constraint set \eqref{eq:israc}, it remains to show that the set \eqref{eq:israc} can be formulated as the intersection of two closed convex sets with known projection operators. In Section \ref{chap:2dim} we demonstrated how to project onto \eqref{eq:israc} for $T=2$ using a geometrical approach. Set \eqref{eq:israc} can be represented as the intersection of $N=T-2$ closed convex sets $\mathcal{U}_k$, where
\begin{subequations}
\begin{gather}\label{eq:setsplit1}
\begin{aligned}
\mathcal{U}_1 \eqdef \bigl\lbrace  u \in \R^{n_u(T\sm 1)}\,\bigl|\,&a^0_{min} \leq u_0 \leq a^0_{max},\\
&\abs{u_1} \leq a,\,\,\abs{u_1 - u_0} \leq r \bigr\rbrace,
\end{aligned}\end{gather}
\vspace{-4pt}
\begin{gather}\label{eq:setsplit2}\begin{aligned}
\mathcal{U}_k \eqdef \bigl\lbrace  u \in \R^{n_u(T\sm 1)}\,\bigl|\,&\abs{u_k} \leq a,\,\,\abs{u_{k-1}} \leq a,\qquad\,\,\\
&\abs{u_k - u_{k-1}} \leq r\bigr\rbrace,
\end{aligned}
\end{gather}\end{subequations}
and $k=2,\dots,T-1$. Let $\pi_k$ denote the procedure presented in Section \ref{chap:2dim} applied element-wise to the $n_u$ elements of $u_{k-1}$ and $u_k$, respectively, with $a_{max}^0\!=\!\sm a_{min}^0\!=\!a$ for $k>1$. Let $\pi_k^1$ and $\pi_k^2$ be the resulting $n_u$ projections of $u_{k-1}$ and $u_k$, respectively. Then the projection operators $\mathcal{P}_{\mathcal{U}_k}$ can be written as
\begin{subequations}\label{eq:projdef}
\begin{align}
\mathcal{P}_{\mathcal{U}_1}(u) &= (\pi_1^1,\pi_1^2,u_3,\dots,u_{T\sm 1})^\Tr,\\
\mathcal{P}_{\mathcal{U}_2}(u) &= (u_0,\pi_2^1,\pi_2^2,u_4,\dots,u_{T\sm 1})^\Tr,\\
&\,\,\,\,\scaleto{\vdots}{10pt}\nonumber\\
\mathcal{P}_{\mathcal{U}_k}(u) &= (u_0,u_1,\dots,\pi_k^1,\pi_k^2,\dots,u_{T\sm 1})^\Tr.
\end{align}\end{subequations}
Let $\mathcal{U}_{e}$ and $\mathcal{U}_{o}$ denote the intersection of sets $\mathcal{U}_i$ grouped by even and odd indices, i.e.\
\begin{subequations}\label{eq:oddeven}\begin{align}
\mathcal{U}_{e} &\eqdef \mathcal{U}_2 \cap \mathcal{U}_4 \cap \dots \cap \mathcal{U}_{T-2},\\
\mathcal{U}_{o} &\eqdef \mathcal{U}_1 \cap \mathcal{U}_3 \cap \dots \cap \mathcal{U}_{T-3},
\end{align}\end{subequations}
where we have assumed for the purposes of explanation that $T$ is even. The projections $\mathcal{P}_{\mathcal{U}_{even}}(u^\circ)$ and $\mathcal{P}_{\mathcal{U}_{odd}}(u^\circ)$ are obtained by combining the corresponding projections from \eqref{eq:projdef} and given by
\begin{subequations}\begin{align}
\mathcal{P}_{\mathcal{U}_{e}}(u) &= (u_0,\pi_2^1,\pi_2^2,\pi_4^1,\pi_4^2,\dots,\pi_{T\sm 1}^1,\pi_{T\sm 1}^2)^\Tr,\\
\mathcal{P}_{\mathcal{U}_{o}}(u) &= (\pi_1^1,\pi_1^2,\pi_3^1(0),\pi_3^2(0),\dots,u_{T\sm 1})^\Tr.
\end{align}\end{subequations}
By setting $\mathcal{P}_{\mathcal{U}_{1}}=\mathcal{P}_{\mathcal{U}_{e}}$ and $\mathcal{P}_{\mathcal{U}_{2}}=\mathcal{P}_{\mathcal{U}_{o}}$, Algorithm \ref{alg:dykstra} can be applied to project onto the input rate and amplitude constraint set \eqref{eq:israc}.

 The choice of using a 2-dimensional projection in combination with Dykstra's method is mainly motivated by the geometrical approach of Section \ref{chap:2dim}. Another possibility would be to compute a formula for a 3-dimensional projection $\tilde{\pi}_k$ onto the set $\mathcal{U}_k\cap\mathcal{U}_{k+1}$ for $n_u=1$ using one of the methods outlined in Section \ref{chap:constraints}. Compared to the 2-dimensional case, fewer $\tilde{\pi}_k$ would have to be evaluated at the expense of increased complexity. How this would affect the computational performance and the convergence rate of the method is not clear a-priori \citep{Han1988}.

\begin{figure}
\begin{center}
\begin{overpic}[scale=0.9]{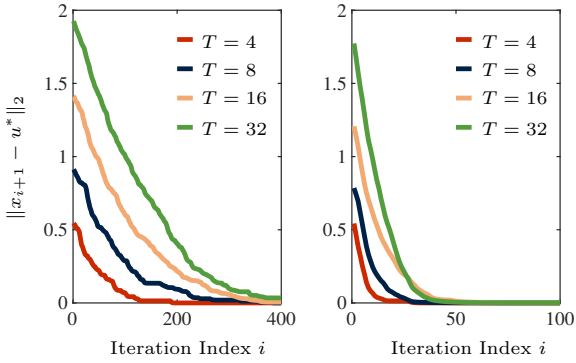}
	% LHS
	\put(32.5,52.5){\scriptsize{$T=4$}}
	\put(32.5,48){\scriptsize{$T=8$}}
	\put(32.5,43.5){\scriptsize{$T=16$}}
	\put(32.5,38.5){\scriptsize{$T=32$}}
	\put(18,3.5){\scriptsize{Iteration Index $i$}}
	\put(1.5,25){\rotatebox{90}{\scriptsize{$\twonorm{x_{i+1}-u^*}$}}}
	% RHS
	\put(77,52.5){\scriptsize{$T=4$}}
	\put(77,48){\scriptsize{$T=8$}}
	\put(77,43.5){\scriptsize{$T=16$}}
	\put(77,38.5){\scriptsize{$T=32$}}
	\put(62.5,3.5){\scriptsize{Iteration Index $i$}}
	%\put(48,25){\rotatebox{90}{\scriptsize{$\twonorm{x_{i+1}-u^*}$}}}
\end{overpic}
\end{center}\caption{Application of Dykstra's method to the set \eqref{eq:israc} with $r\!=\!a\!=\!1$. 100 starting points $u^\circ$ were drawn from a random normal distribution with zero mean and a standard deviation of 100 (left) and 10 (right), respectively. Shown is the distance of variable $x_{i+1}$ from Algorithm \ref{alg:dykstra} to an accurate projection $u^*$.}\label{fig:projerror}
\end{figure}

%\subsection{Numerical Example}\label{chap:dykstrasim}
 Figure \ref{fig:projerror} compares the output of Algorithm \ref{alg:dykstra} applied to the input rate and amplitude constraint set \eqref{eq:israc} with $r=a=1$ for different horizons $T$. Shown is the distance of iterates $x_{i+1}$ from Algorithm \ref{alg:dykstra} to the solution $u^*$ obtained using an interior-point method. For the figure we selected 100 starting points $u^\circ$ from a normal distribution with zero mean and a standard deviation of $100$ (left) and $10$ (right), respectively.

 %The simulations in Figure \ref{fig:projerror} suggest that the required number of iterations is proportional to the initial distance $\twonorm{u^\circ - u^*}$ of the starting point $u^\circ$ to its projection $u^*$. Since Algorithm \ref{alg:dykstra} will be used in conjunction with the FGM (Algorithm \ref{alg:fgm1}), the initial distance $\twonorm{u^\circ - u^*}$ will decrease with the iterates of the FGM. Fewer iterates will be required to solve the projection as the FGM advances.

\section{MPC for Systems with Input Rate and Amplitude Constraints}\label{chap:mpcfgm}
%\subsection{Algorithm}\label{chap:mpcfgm}
The FGM (Algorithm \ref{alg:fgm1}) and Dykstra's projection method (Algorithm \ref{alg:dykstra}) are combined in Algorithm \ref{alg:fgm2}, where the projection $\mathcal{P}_{\mathcal{U}}$ onto the input rate and amplitude constraint set was replaced with Dykstra's algorithm. Compared to the ADMM implementation, Algorithm \ref{alg:fgm2} bears several advantages. Dykstra's algorithm makes the augmentation of decision variables superfluous. While the ADMM formulation \eqref{eq:mpcqpaug} requires $2T-1$ decision variables, Algorithm \ref{alg:fgm2} reduces the number of decision variables to $T$. This not only greatly reduces the computation time, as the next section will show, but also lowers the memory footprint. If we assume that all matrices are dense and neglect the storage of vectors, then Algorithm \ref{alg:fgm2} reduces the memory footprint by approximately $T^2 / (T^2+(T-1)^2 n_v / n_u)$, which roughly amounts to a reduction of $50\%$ for set \eqref{eq:israc} and large $T$. Dykstra's algorithm barely introduces any memory footprint because it solely involves Boolean operations and vector additions. Moreover, arrays allocated by Algorithm \ref{alg:fgm1} can be used as temporary placeholders to execute Dykstra's method. Stand-alone C-language implementations of Algorithm \ref{alg:fgm1} and \ref{alg:dykstra} can be found under \citep{BRPGITHUB, FGMGITHUB}. To this end, note that Algorithms \ref{alg:admm}, \ref{alg:fgm1} and \ref{alg:fgm2} could be warm-started using the solution computed at time step $t-1$, e.g.\ setting $y_1 = u_{\sm 1} = u(t-1)$ in Algorithm \ref{alg:fgm2}.

\begin{algorithm}[H]
 \caption{Fast Gradient Method and Dykstra's Algorithm applied to MPC problem \eqref{eq:mpc} with set \eqref{eq:israc}}\label{alg:fgm2}
 \begin{algorithmic}[1]
 \Require{$x_0=x(t)$ and $u_{\sm 1}=u(t-1)$}
 \Ensure{$u(t)=u_{\scaleto{I}{4pt}_{max}}$}
  \State{Set $q=G^\Tr H x_0$ and $y_1=u_1=0$}
  \State{Initialize projection using $u_{\sm 1}$}
  \For{$i = 1$ to $I_{max}$}
  	\State $t_i = (I - J/\lambda_{max})y_i - q / \lambda_{max}$
  	\State Set $x_0=t_i$, $\mu_0=0$ and $\gamma_0=0$
  	\For{$j = 1$ to $J_{max}$}
      \State{$w_j=\mathcal{P}_{\mathcal{U}_{e}}(x_j+\mu_j)$}
      \State{$\mu_{j+1}=x_j+\mu_j-w_j$}
      \State{$x_{j+1}=\mathcal{P}_{\mathcal{U}_{o}}(w_j+\gamma_j)$}
      \State{$\gamma_{j+1}=w_j+\gamma_j-w_{j+1}$}
      \If {$\infnorm{x_{j+1}-x_{j}} < \epsilon$}
  		\Break
  	  \EndIf
  	\EndFor
  	\State $u_{i+1} = x_{j+1}$
  	\State $y_{i+1} = (1+\beta)u_{i+1} - \beta u_i$
  \EndFor
 \end{algorithmic}
\end{algorithm}

\subsection{Numerical Studies}\label{chap:mpcfgmsim}
Algorithm \ref{alg:fgm2} is compared against an ADMM implementation which -- as described in Section \ref{chap:admm} -- commonly uses an augmentation of decision variables to simplify the projection. While Algorithm \ref{alg:fgm2} was implemented in C and can be found under \citet{BRPGITHUB, FGMGITHUB}, the Operator Splitting Quadratic Program solver \citep{OSQPMAIN} -- a constrained QP solver that uses ADMM -- was employed for Algorithm \ref{alg:admm}. In order to avoid refactoring the matrix on the left-hand side of step 3 of Algorithm \ref{alg:admm}, our benchmark ADMM implementation uses a constant penalty parameter $\rho$. Neither of the C-programs includes a non-standard C-library.

\begin{figure}
\begin{center}
\begin{overpic}[scale=0.9]{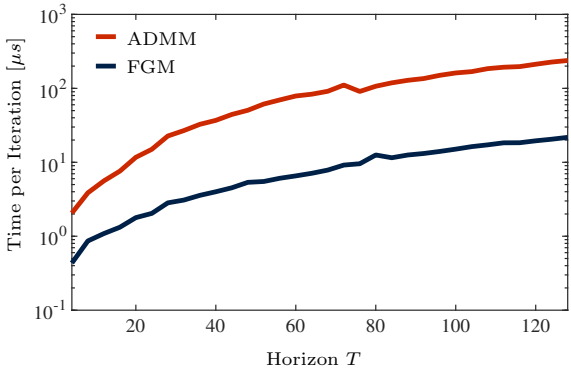}
	\put(22,54){\scriptsize{ADMM}}
	\put(22,49.5){\scriptsize{FGM}}

	\put(44,3){\scriptsize{Horizon $T$}}
	\put(3,21){\rotatebox{90}{\scriptsize{Time per Iteration $\left[\mu s\right]$}}}
\end{overpic}
\end{center}\caption{Comparison of the logarithmically-scaled average execution times of one iteration of Algorithm \ref{alg:admm} (ADMM, red) and Algorithm \ref{alg:fgm2} (FGM, blue), respectively, applied to problem \eqref{eq:mpcqp} with set \eqref{eq:israc}. Both algorithms were benchmarked using a stand-alone C-language implementation and using 100 randomly generated problems per horizon $T$. The FGM includes Dykstra's method with a termination check on every \nth{10} iterate.}\label{fig:avgtimeperiter}
\end{figure}

 Figure \ref{fig:avgtimeperiter} compares the average execution times of one iteration of the C-language implementation of the ADMM and the FGM applied to problem \eqref{eq:mpcqp} with the input rate and amplitude constraint set $\mathcal{U}$ as defined in \eqref{eq:israc} with $r=a=1$. The OSQP solver uses a matrix-factorization to solve step 3 of Algorithm \ref{alg:admm} and the time for factorizing the matrix is excluded from Figure \ref{fig:avgtimeperiter}. The problem data $(J,q)$ is randomly generated and the average execution times are benchmarked over $100$ problems per horizon. Figure \ref{fig:avgtimeperiter} reveals that the combination of the FGM and Dykstra's method greatly reduces the computation time for one solver iteration. The performance gain is due to the fact that Dykstra's method makes the augmentation of decision variables \eqref{eq:mpcqpaug} unnecessary, which in case of the input rate and amplitude set amounts to tripling the number of decision variables. Because Dykstra's method involves only vector additions and Boolean operations, the projection algorithm only requires a few processor cycles. The termination criterion for Dykstra's algorithm is checked on every \nth{10} iterate.

\begin{figure}
\begin{center}
\begin{overpic}[scale=0.9]{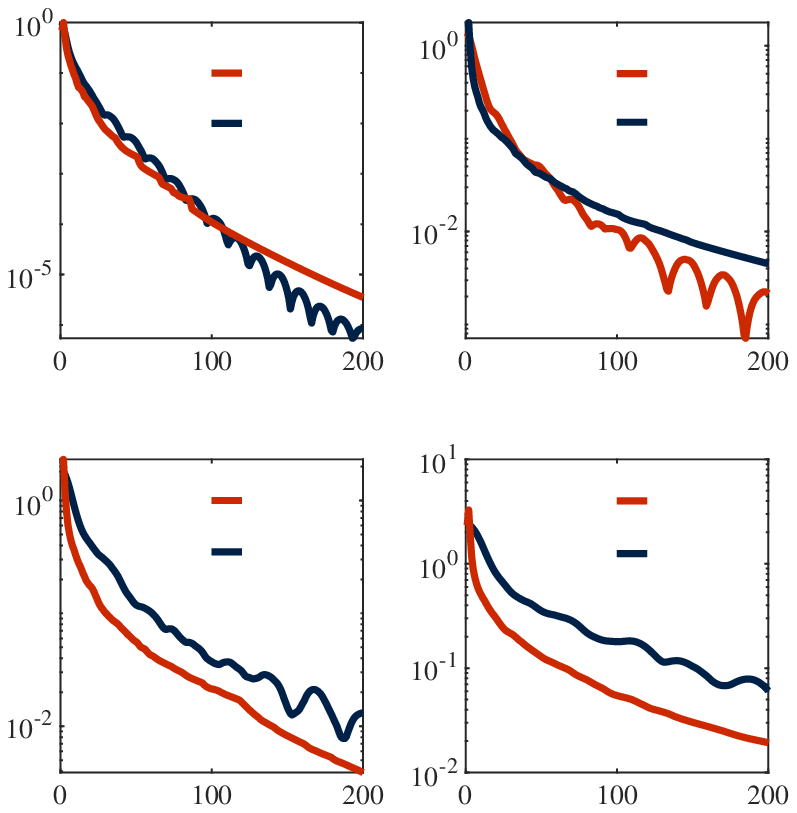}
	\put(18,60){\normalsize{$T=4$}}
	\put(59.5,60){\normalsize{$T=8$}}
	\put(18,15){\normalsize{$T=16$}}
	\put(59.5,15){\normalsize{$T=32$}}

	\put(33.5,83.5){\scriptsize{ADMM}}
	\put(33.5,78){\scriptsize{FGM}}
	\put(75,83.5){\scriptsize{ADMM}}
	\put(75,78){\scriptsize{FGM}}

	\put(33.5,39.5){\scriptsize{ADMM}}
	\put(33.5,34){\scriptsize{FGM}}
	\put(75,39.5){\scriptsize{ADMM}}
	\put(75,34){\scriptsize{FGM}}

	\put(19.5,5){\scriptsize{Iteration Index $i$}}
	\put(61.5,5){\scriptsize{Iteration Index $i$}}
	\put(19.5,49){\scriptsize{Iteration Index $i$}}
	\put(61.5,49){\scriptsize{Iteration Index $i$}}

	\put(4,22){\rotatebox{90}{\scriptsize{$\twonorm{u_i - u^*}$}}}
	\put(4,67){\rotatebox{90}{\scriptsize{$\twonorm{u_i - u^*}$}}}
\end{overpic}
\end{center}\caption{Practical convergence behaviour of Algorithm \ref{alg:admm} (ADMM, red) and Algorithm \ref{alg:fgm2} (FGM, blue) applied to problem \eqref{eq:mpcqp} with set \eqref{eq:israc} for $T=\lbrace 4, 8, 16, 32\rbrace$. Shown is the distance between the iterates and a solution $u^*$ computed using an interior-point method. Algorithm \ref{alg:admm} uses a constant penalty parameter $\rho$. The problem data is randomly generated and the distances are averaged over $100$ problems.}\label{fig:avgconvergence}
\end{figure}

 Figure \ref{fig:avgconvergence} compares the practical convergence behaviour of Algorithm \ref{alg:admm} (ADMM, red) and Algorithm \ref{alg:fgm2} (FGM, blue) applied to problem \eqref{eq:mpcqp} with the input rate and amplitude constraint set as defined in \eqref{eq:israc} with $r=a=1$. Depicted is the average distance between a high-accuracy solution $u^*$ calculated using an interior-point method and the solution at iteration $i$ of Algorithm \ref{alg:admm} and \ref{alg:fgm2}, respectively. The problem data $(J,q)$ is randomly generated and the distances are averaged over $100$ problems per horizon $T$. For Algorithm \ref{alg:fgm2}, Dykstra's method is implemented as follows: An initial verification is applied to avoid executing the algorithm for vectors $t_i$ that lie inside set \eqref{eq:israc}. If $t_i\notin\mathcal{U}$, Dykstra's method is run for a fixed number of iterations $J_{max}=50$.

\subsection{Example}\label{chap:example}
We study the application of Algorithms \ref{alg:admm} and \ref{alg:fgm2} to an aircraft stabilization problem. The linearized model is taken from \citep{AIRCRAFTMODEL}, discretized with a sampling time of $T_s =1ms$ and an MPC problem formulated in its condensed form \eqref{eq:mpcqp} for horizons $T=\lbrace 4, 8, 16, 32\rbrace$ and with $Q=R=I$. The aircraft model features $n_u=3$ inputs, i.e.\ $u_k = (u_k^{ail}, u_k^{stab}, u_k^{rud})^\Tr$, where $u_k^{ail}, u_k^{stab}$ and $u_k^{rud}$ denote the differential aileron deflection, the differential stabilizer deflection and the rudder deflection, respectively. The inputs are constrained by rate and amplitude constraints with $(a^{ail}, a^{stab}, a^{rud}) = (25, 24, 30)$ and $(r^{ail}, r^{stab}, r^{rud}) = (200 T_s, 80 T_s, 82 T_s)$, respectively. 

While Algorithms \ref{alg:admm} and \ref{alg:fgm2} are used to solve the condensed MPC problem \eqref{eq:mpcqp}, the dynamics of the aircraft are simulated using \eqref{eq:mpcB} starting with a random non-zero initial condition. The closed-loop simulation is run for 1000 time steps. In order to be able to compare the solution accuracies, a termination criterion for Algorithm \ref{alg:fgm2} is introduced. In addition, it is verified that both algorithms produce a similar closed-loop behavior.

\begin{figure}
\begin{center}
\begin{overpic}[scale=0.9]{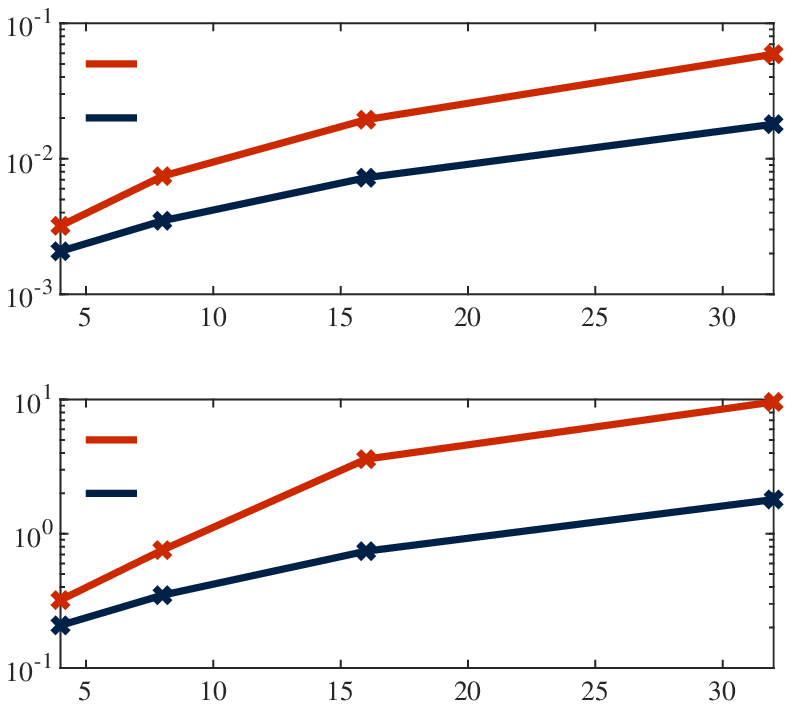}
	\put(23,75.7){\scriptsize{ADMM}}
	\put(23,70){\scriptsize{FGM}}
	
	\put(23,35.7){\scriptsize{ADMM}}
	\put(23,30){\scriptsize{FGM}}

	\put(44,5){\scriptsize{Horizon $T$}}
	\put(44,45){\scriptsize{Horizon $T$}}
	
	\put(3,49){\rotatebox{90}{\scriptsize{Time per Iteration $\left[ms\right]$}}}
	\put(3,15){\rotatebox{90}{\scriptsize{Total Time $\left[ms\right]$}}}
\end{overpic}
\end{center}\caption{Comparison of the logarithmically-scaled execution times of Algorithm \ref{alg:admm} (ADMM, red) and Algorithm \ref{alg:fgm2} (FGM, blue), respectively, applied to the aircraft stabilization problem of Section \ref{chap:example}. The top and bottom row show the average execution time per iteration and the total execution time to solve the MPC problem \eqref{eq:mpcqp}, respectively.}\label{fig:example}
\end{figure}

Figure \ref{fig:example} depicts the average time per iteration (top) and the average total time (bottom) required to solve one instance of problem \eqref{eq:mpc}. Both rows are averaged over the 1000 time steps of the closed-loop simulation. Compared to Figure \ref{fig:avgtimeperiter} where $n_u=1$ was used, it can be seen that the $n_u=3$ inputs introduce an overhead and slow down Algorithm \ref{alg:fgm2} for smaller horizons. The performance advantage is approximately regained for $T=32$. The bottom row confirms the results from Figure \ref{fig:avgconvergence}, which showed that a similar convergence behavior can be expected from Algorithms \ref{alg:admm} and \ref{alg:fgm2}.

%%%%%%%%%%%%%%%%%%%%%%%%%%%%%%%%%%%%%%%%%%%%%%%%%%%%%%%%%%%%%%%
% SECTION 6: CONCLUSIONS
%%%%%%%%%%%%%%%%%%%%%%%%%%%%%%%%%%%%%%%%%%%%%%%%%%%%%%%%%%%%%%%
\section{Conclusions and Future Work}
This paper demonstrated the use of Dykstra's method to project onto the input rate and amplitude constraint set. A procedure to solve the 2-dimensional projection was presented and a C-language implementation provided. It was shown how the input constraint set can be represented as the intersection of two sets. Using the formula for the 2-dimensional case, the projection was extended to higher dimensions and computed iteratively using Dykstra's method. Several simulations showed the fast convergence of the projection algorithm. The choice of using an underlying 2-dimensional projection was motivated by a geometrical approach. It is certain that choosing a different dimension for the underlying projection would affect the performance of the algorithm. Investigating the performance for different underlying projections is subject of future work.

 In order to solve a linear MPC problem constrained by input rate and amplitude limits, Dykstra's method was embedded in an FGM that uses a constant step scheme. The resulting algorithm was compared against an ADMM implementation. The ADMM uses a decision variable augmentation to accommodate the constraints. Using a C-language implementation applied to several example problems, it was shown that the combination of the FGM and Dykstra's algorithm significantly reduces the computation time compared to an ADMM implementation. Moreover, the practical convergence behaviour was simulated and it was shown that the combination of FGM and Dykstra's method converges as quickly as the ADMM implementation. In addition, it was demonstrated that the combined algorithms approximately halve the memory footprint, which is of particular importance for an implementation on an embedded system.

 While it has been shown that the practical convergence behaviour of Algorithm \ref{alg:fgm2} matches the one of Algorithm \ref{alg:admm}, future work will focus on establishing convergence criteria for the combination of Algorithms \ref{alg:fgm1} and \ref{alg:dykstra}. Dykstra's method introduces a projection error that can be made arbitrarily small by increasing the number of iterations \citep[Chapter 29.1]{BAUSCHKEBOOK}. The convergence of ADMM under an inexact projection has been proved in \citet[Theorem 8]{Eckstein1992}, where it was shown that the algorithm converges provided that the sum of the projection errors remains bounded. Since the projection error of Dykstra's algorithm can be made arbitrarily small, the boundedness of the sum of the projection errors can be enforced and it might be expected that a similar proof can be formulated for Algorithm \ref{alg:fgm2}. A practically relevant situation is when the FGM as well as Dykstra's method are run for a fixed number of iterations. A proof exists for the FGM with an exact projection in \citet{CERTIFMPC}, where it was shown that a certain level of suboptimality can be guaranteed provided that the initial condition $x_0$ of the MPC problem \eqref{eq:mpc} lies within a bounded set. Future efforts will also aim at extending the proof from \citet{CERTIFMPC} to the case of a projection using Dykstra's method.

%%%%%%%%%%%%%%%%%%%%%%%%%%%%%%%%%%%%%%%%%%%%%%%%%%%%%%%%%%%%%%%
% BIBLIOGRAPHY
%%%%%%%%%%%%%%%%%%%%%%%%%%%%%%%%%%%%%%%%%%%%%%%%%%%%%%%%%%%%%%%
%\bibliographystyle{IEEEtran}
\small
\bibliography{bibliography}

\end{document}